\documentclass{article}

\usepackage{amsmath,amssymb,amsfonts,theorem}
\usepackage{graphicx,color}
\usepackage{subfigure}
\graphicspath{{./Figures/}}

{ \theorembodyfont{\normalfont} 

}

\newtheorem{theorem}{Theorem}
\newtheorem{lemma}{Lemma}

\newtheorem{proposition}{Proposition}

\def\scr#1{{\cal #1}}

\newcommand{\bigintersect}{\bigcap}

\newcommand{\reals}{\mathbb{R}} 

\newcommand{\R}{\reals}
\newcommand{\dfb}{\stackrel{\Delta}{=}}

\newcommand{\dst}{\displaystyle}

\def\be{\begin{equation}}
\def\ee{\end{equation}}
\def\ba{\begin{array}}
\def\ea{\end{array}}
\def\eqa{\begin{eqnarray}}
\def\eqe{\end{eqnarray}}

\author{ \parbox{4 in}{\centering Claudio De Persis,
          Hui Liu, and Ming Cao
          \thanks{C. De Persis is with the Faculty of Engineering Technology (CTW), University of Twente, the Netherlands ({\tt\small
         c.depersis@utwente.nl}) and with the Department of Computer and System Sciences, Sapienza University of Rome, Italy. H. Liu and M. Cao are with the Faculty of Mathematics and Natural
Sciences, ITM, University of Groningen, the Netherlands
({\tt\small \{hui.liu, m.cao\}@rug.nl}).}
          } }

\title{\LARGE \bf
Control of one-dimensional guided formations\\ using coarse
information}

\begin{document}

\maketitle

\begin{abstract}
Motivated  by applications in intelligent highway systems, the
paper studies the problem of guiding mobile agents in a
one-dimensional formation to their desired relative positions.
Only coarse information is used which is communicated from a
guidance system that monitors in real time the agents' motions.
The desired relative positions are defined by the given distance
constraints between the agents under which the overall formation
is rigid in shape and thus admits locally a unique realization. It
is shown  that  even when the guidance system can only transmit at
most four bits of information to each agent, it is still possible
to design control laws to guide the agents to their desired
positions. We further delineate the thin set of initial conditions
for which the proposed control law may fail using the example of a
three-agent formation.
Tools from non-smooth analysis are utilized  for the convergence
analysis.

\end{abstract}

\section{Introduction} \label{se:intro}
In recent years, various  ideas have been proposed to realize
intelligent highway systems to reduce traffic congestions and
improve safety levels. It is envisioned that navigation,
communication and automatic driver assistance systems are critical
components \cite{BrVa97,Bi05,Sh07}. A great deal of monitoring and
controlling capabilities have been implemented through  roadside
infrastructures, such as cameras, sensors, and control and
communication stations. Such systems can work together to monitor
in real time the  situations on  highways and at the same time
guide vehicles to move in a coordinated fashion, e.g. to keep
appropriate distances from the vehicles in front of and behind
each individual vehicle. In  intelligent highway systems, the
guiding commands are expected to be simple and formatted as short
digital messages  to scale with the number of vehicles and also to
avoid conflict with the automatic driver assistance systems
installed within the vehicles. Similar guided formation control
problems also arise when navigating mobile robots or docking
autonomous vehicles \cite{TuCrTrMoBoDo06}.

Motivated by this problem  of guiding platoons of vehicles on
highways, we study in this paper the problem of controlling a
one-dimensional multi-agent formation using only \emph{coarsely}
quantized information. The formation to be considered are rigid
under inter-agent distance constraints and thus its shape is
uniquely determined locally. Most of the existing work on
controlling rigid formations of mobile agents, e.g.
\cite{AnYuFiHe08,CaMoCDC07,KrBrFr08}, assumes that there is no
communication bandwidth constraints and thus real-valued control
signals are utilized. The idea of quantized control through
digital communication channels has been applied to consensus
problems, e.g.\ \cite{KaBaSr07,FrCaFaZa09} and references therein,
and more recently to formation control problems \cite{DiJo10}. The
uniform quantizer and logarithmic quantizer \cite{NaFaZaEv07} are
among the most popular choices for designing such controllers with
quantized information. Moreover, the paper \cite{CePeFr10} has
discussed Krasowskii solutions and hysteretic quantizers in
connection with continuous-time average consensus algorithms under
quantized measurements.

The problem studied in this  paper distinguishes itself from the
existing work in that it explores the limit of the least bandwidth
for controlling a one-dimensional rigid formation by using a
quantizer in its simplest form with  only two quantization levels.
As a result, for each agent in the rigid formation, at most four
bits of bandwidth is needed for the communication with the
navigation controller. The corresponding continuous-time model
describing the behavior of the overall multi-agent formation is,
however, non-smooth and thus an appropriate notion of solution
\cite{CePe07} has to be defined first. We use both the Lyapunov
approach and trajectory-based approach to prove convergence since
the former provides a succinct view about the dynamic behavior
while the latter leads to insight into the set of initial
positions for which the proposed controller may fail. We also
discuss some situations when different assumptions about the
quantization scheme are made and indicate those scenarios in which
the formation control problem with quantized information can be
challenging to solve.

The rest of the paper is  organized as follows. We first formulate
the one-dimensional guided formation control problem with coarsely
quantized information in section \ref{se:formulation}. Then in
section \ref{se:analysis}, we provide the convergence analysis
results first using the Lyapunov method and then the
trajectory-based method.
Simulation results are presented in section \ref{se:simulation} to
validate the theoretical analysis. We make concluding remarks in
section \ref{se:conclusion}.

\section{Problem formulation} \label{se:formulation}
The one-dimensional guided formation that we are interested in
consists of $n$ mobile agents. We consider the case when the
formation is rigid \cite{AnYuFiHe08}; to be more specific,  if we
align the given one-dimensional space with the $x$-axis in the
plane and label the agents along the positive direction of the
$x$-axis by $1, \ldots, n$, then the geometric shape of the
formation is  specified by the given pairwise distance constraints
$|x_i-x_{i+1}| = d_i$, $i = 1, \ldots, n-1$, where $d_i>0$ are
desired distances. Although the guidance system can monitor the
motion of the agents in real time, we require that it can only
broadcast to the mobile agents quantized guidance information
through digital channels. In fact, we explore the limit for the
bit constraint by utilizing the quantizer that only has two
quantization levels and consequently its output only takes up one
bit of bandwidth. The quantizer that is under consideration takes
the form of the following sign function: For any $z\in \R$,
\[
\textrm{sgn}(z)=\left\{ \ba{ll}
+1 & z\ge 0\\
-1 & z<0\;. \ea\right.
\]
Each agent, modeled by a kinematic point, then moves according to
the following rules utilizing the coarsely quantized information:
\be\label{quantized.n.agent.system} \ba{rcl}
\dot x_1 &=& -k_1 \textrm{sgn} (x_1-x_2) \textrm{sgn} (|x_1-x_2|-d_1)\\
\dot x_i &=& \textrm{sgn}(x_{i-1}-x_i) \textrm{sgn}(|x_{i-1}-x_i|-d_{i-1})-\\
&& k_i \textrm{sgn}(x_i-x_{i+1}) \textrm{sgn}(|x_i-x_{i+1}|-d_{i}),\\ &&\hfill i=2,\ldots,n-1\\
\dot x_n &=& \textrm{sgn}(x_{n-1}-x_n)
\textrm{sgn}(|x_{n-1}-x_n|-d_{n-1}) \ea \ee
where $x_i\in \R$ is the position of agent $i$ in the
one-dimensional space aligned with the $x$-axis, and $k_i>0$ are
gains to be designed. Note that since each agent is governed by at
most two distance constraints, as is clear from
(\ref{quantized.n.agent.system}), a bandwidth of four bits is
sufficient for the communication between the guidance system and
the agents $2, \ldots, n-1$ and the required bandwidths for the
guidance signals for agents $1$ and $n$ are both 2 bits. Hence, in
total only $4n-2$ bits of bandwidth is used.

The main goal of this paper is to demonstrate under this extreme
situation of using coarsely quantized information, the formation
still exhibits satisfying convergence properties under the
proposed maneuvering rules. Towards this end, we introduce the
variables of relative positions among the agents \be\label{z}
z_i\dfb x_i-x_{i+1}\;,\quad i=1,2,\ldots,n-1\;. \ee Let us express
the system in the $z$-coordinates to obtain
\be\label{quantized.n.agent.system.z} \ba{rcl}
\dot z_1 &=& -(k_1+1) \textrm{sgn}(z_1) \textrm{sgn}(|z_1|-d_1)\\[2mm] && + k_2 \textrm{sgn}(z_2) \textrm{sgn}(|z_2|-d_2)\\[2mm]
\dot z_i &=& \textrm{sgn}(z_{i-1}) \textrm{sgn}(|z_{i-1}|-d_{i-1}) \\[2mm]&& -(k_i+1) \textrm{sgn}(z_i) \textrm{sgn}(|z_i|-d_{i})\\[2mm]
&& +k_{i+1} \textrm{sgn}(z_{i+1}) \textrm{sgn}(|z_{i+1}|-d_{i+1}),\\[2mm] &&   \hfill i=2,\ldots,n-2\\[2mm]
\dot z_{n-1} &=& \textrm{sgn}(z_{n-2})
\textrm{sgn}(|z_{n-2}|-d_{n-2})\\ [2mm] &&-(k_{n-1}+1)
\textrm{sgn}(z_{n-1}) \textrm{sgn}(|z_{n-1}|-d_{n-1})\;. \ea \ee
To study the dynamics of the system above, we need to first
specify what we mean by the solutions of the system. Since the
vector field $f(z)$ on the right-hand side is discontinuous, we
consider Krasowskii solutions, namely solutions to the
differential inclusion $\dot z\in {\cal K}(f(z))$, where
\[
{\cal K}(f(z))=\bigintersect_{\delta
>0}{\overline {\rm co}}\, (f (B(z,\delta)))\;,
\]
$\overline {\rm co}$ denotes the involutive closure of a set, and
$B(z,\delta)$ is the ball centered at $z$ and of the radius
$\delta$. The need to consider these solutions becomes evident in
the analysis in the next section. Since the right-hand side of
(\ref{quantized.n.agent.system}) is also discontinuous, its
solutions are to be intended in the Krasowskii sense as well. Then
we can infer conclusions on the behavior of
(\ref{quantized.n.agent.system}) provided that each solution $x$
of (\ref{quantized.n.agent.system}) is such that $z$ defined in
(\ref{z}) is a Krasowskii solution of
(\ref{quantized.n.agent.system.z}). This is actually the case by
\cite{BP-SSS:87}, Theorem 1, point 5), and it is the condition
under which we consider (\ref{quantized.n.agent.system.z}). It
turns out  that the $z$-system (\ref{quantized.n.agent.system.z})
is easier to work with for the convergence analysis that we
present in detail in the next section.

\section{Convergence analysis} \label{se:analysis}
In this section, after identifying the equilibria of the system,
we present two different approaches for convergence analysis. The
first is based on a Lyapunov-like function and the second examines
the vector field in the neighborhood of the system's trajectories.

\subsection{Equilibria of the system}
We start the analysis of system (\ref{quantized.n.agent.system.z})
by looking at the discontinuity points of the system.  A
discontinuity point is a point at which the vector field on the
right-hand side of the equations above is discontinuous. Hence,
the set ${\cal D}$ of all the discontinuity points is:
\[
{\cal D}=\{z\in \R^{n-1}: \Pi_{i=1}^{n-1} z_i (|z_i|-d_{i})=0\}\;.
\]
It is  of interest to characterize the set of equilibria:
\begin{proposition}\label{lemma.equilibria}
Let $k_1+1>k_2$, $k_i>k_{i+1}$  for $i=2,\ldots,n-2$, and
$k_{n-1}>0$.  The set of equilibria, i.e.\ the set of points for
which $\mathbf{0}\in {\cal K}(f(z))$ with $f(z)$ being the vector
field on the right-hand side of
(\ref{quantized.n.agent.system.z}), is given by
\[
{\cal E}=\{z\in \R^{n-1}: \sum_{i=1}^{n-1} |z_i|
||z_i|-d_{i}|=0\}\;.
\]
\end{proposition}

The proof of this proposition relies on the following lemma.

\begin{lemma} \label{lemma.claim}
For $i\in \{2,\ldots, n-2\}$, if $|z_j|\,||z_j|-d_{j}|=0$ for
$j=1,2,\ldots, i-1$,  and $\mathbf{0}\in {\cal K}(f(z))$, then
$|z_{i}|\, ||z_{i}|-d_{i}|=0$.
\end{lemma}

\medskip

\textit{Proof:} Suppose by contradiction that $|z_{i}|\,
||z_{i}|-d_{i}|\ne 0$. Observe that $z$ belongs to a discontinuity
surface where in particular $|z_{i-1}|\, ||z_{i-1}|-d_{i-1}|=0$.
This implies that in a neighborhood of this point, the state space
is partitioned into different regions where $f(z)$ is equal to
constant vectors. In view of (\ref{quantized.n.agent.system.z}),
the component $i$  of these vectors is equal to one of the
following values: $1-(k_i+1)+k_{i+1}$, $1-(k_i+1)-k_{i+1}$,
$-1-(k_i+1)+k_{i+1}$, $-1-(k_i+1)-k_{i+1}$, if $\textrm{sgn}(z_i)
\textrm{sgn}(|z_i|-d_{i})=1$, or $1+(k_i+1)+k_{i+1}$,
$1+(k_i+1)-k_{i+1}$, $-1+(k_i+1)+k_{i+1}$, $-1+(k_i+1)-k_{i+1}$,
if $\textrm{sgn}(z_i) \textrm{sgn}(|z_i|-d_{i})=-1$. Any $v\in
{\cal K}(f(z))$ is such that its component $i$ belongs to (a
subinterval of) the interval
$[-1-(k_i+1)-k_{i+1},1-(k_i+1)+k_{i+1}]$ if $\textrm{sgn}(z_i)
\textrm{sgn}(|z_i|-d_{i})=1$ (respectively, to the interval
$[-1+(k_i+1)-k_{i+1},1+(k_i+1)+k_{i+1}]$ if $\textrm{sgn}(z_i)
\textrm{sgn}(|z_i|-d_{i})=-1$). In both cases, if $k_i>k_{i+1}$,
then the interval does not contain $0$ and this is a
contradiction. This ends the proof of the lemma. \hfill $\square$

\textit{Proof of Proposition \ref{lemma.equilibria}:} First we
show that if $\mathbf{0}\in {\cal K}(f(z))$, then $z\in {\cal E}$.
As a first step, we observe that $\mathbf{0}\in {\cal K}(f(z))$
implies $|z_1| ||z_1|-d_{1}|=0$. In fact, suppose by contradiction
that the latter is not true. This implies that at the point $z$
for which $\mathbf{0}\in {\cal K}(f(z))$, any $v\in {\cal
K}(f(z))$ is such that the first component takes values in the
interval $[-(k_1+1)-k_2, -(k_1+1)+k_2]$, or in the interval
$[(k_1+1)-k_2, (k_1+1)+k_2]$. In both cases, if $k_1+1>k_2$, then
$0$ does not belong to the interval and this contradicts that
$\mathbf{0}\in {\cal K}(f(z))$. Hence, $|z_1| ||z_1|-d_{1}|=0$.

%

This and Lemma \ref{lemma.claim}  show that  $|z_{i}|\,
||z_{i}|-d_{i}|=0$ for $i=1,2,\ldots, n-2$. To prove that also
$|z_{n-1}|\, ||z_{n-1}|-d_{n-1}|=0$, consider the last equation of
(\ref{quantized.n.agent.system.z}), and again suppose by
contradiction that $|z_{n-1}|\, ||z_{n-1}|-d_{n-1}|\ne 0$. Then
the last component of $v\in {\cal K}(f(z))$ belongs to a
subinterval of $[-1-(k_{n-1}+1), 1-(k_{n-1}+1)]$ or to a
subinterval of $[-1+(k_{n-1}+1), 1+(k_{n-1}+1)]$. If $k_{n-1}>0$,
then neither of these intervals contain $0$ and this is again a
contradiction. This concludes the first part of the proof, namely
that $\mathbf{0}\in {\cal K}(f(z))$ implies $z\in {\cal E}$.

Now  we let $z\in {\cal E}$ and prove that $\mathbf{0}\in {\cal
K}(f(z))$. By definition, if $z\in {\cal E}$, then $z$ lies at the
intersection of $n-1$ planes, which partition $\R^n$ into
$\nu\stackrel{\cdot}{=} 2^{n-1}$ regions, on each one of which
$f(z)$ is equal to a different constant vector. Any $v \in {\cal
K}(f(z))$ is the convex combination of these $\nu$ vectors, which
we call $v^{(1)}, \ldots, v^{(\nu)}$. We construct $v \in {\cal
K}(f(z))$ such that $v=\mathbf{0}$. We observe first that, the
component $1$ of the vectors $v^{(i)}$'s can take on four possible
values, namely $(k_1+1)+k_2$, $(k_1+1)-k_2$, $-(k_1+1)+k_2$,
$-(k_1+1)-k_2$, and that there are exactly $\frac{\nu}{4}$ (we are
assuming that $n\ge 3$, as the case $n=2$ is simpler and we omit
the details) vectors among $v^{(1)},\ldots,v^{(\nu)}$ whose first
component is equal to $(k_1+1)+k_2$, $\frac{\nu}{4}$ whose first
component is equal to $(k_1+1)-k_2$ and so on. As a consequence,
if $\lambda_i=\frac{1}{\nu}$ for all $i=1,2,\ldots,\nu$, then
$\sum_{j=1}^\nu \lambda_j v_1^{(j)}=0$.\\
Similarly, the component $i$, with $i=2,\ldots, n-2$, can take on
eight possible values ($1+(k_i+1)+k_{i+1}$,
$1+(k_i+1)-k_{i+1}$,$\ldots, -1-(k_i+1)-k_{i+1}$ -- see the
expression of $\dot z_i$ in (\ref{quantized.n.agent.system.z}))
and as before, the set $v^{(1)},\ldots,v^{(\nu)}$ can be
partitioned into $\frac{\nu}{8}$ sets, and each vector in a set
has the component $i$ equal to one and only one of the eight
possible values. Moreover, these values  are such that
$\sum_{j=1}^\nu
\lambda_j v_i^{(j)}=0$. \\
Finally, if $i=n-1$, the set $v^{(1)},\ldots,v^{(\nu)}$ can be
partitioned into four sets, and each vector in a set has the last
component equal to one and only one of the four possible values
$1+(k_{n-1}+1)$, $1-(k_{n-1}+1)$, $-1+(k_{n-1}+1)$,
$-1-(k_{n-1}+1)$. Hence, $\sum_{j=1}^\nu \lambda_j
v_{n-1}^{(j)}=0$. Let now $v\in {\cal K}(f(z))$ be such that
$v=\sum_{i=1}^{\nu} \lambda_i v^{(i)}$, with
$\lambda_i=\frac{1}{\nu}$ for all $i$. Since $\sum_{j=1}^\nu
\lambda_j v_i^{(j)}=0$ for all $i=1,2,\ldots, n-1$, then
$v=\mathbf{0}$ and this proves that for all $z\in {\cal E}$, we
have $\mathbf{0}\in {\cal K}(f(z))$. This completes the proof.
\hfill $\square$

Next, we show that the equilibrium set $\scr{E}$ is attractive.

\subsection{Lyapunov function based analysis}
Now we are in a position to present the main convergence result.
\begin{theorem}
If \be\label{gains} k_1\ge k_2 \;,\; k_i\ge k_{i+1}+1\;,\;
i=2,\ldots,n-2\;,\; k_{n-1}\ge 1\;, \ee then all the Krasowskii
solutions to (\ref{quantized.n.agent.system.z}) converge to (a
subset  of) the equilibria set ${\cal E}$.
\end{theorem}

\textit{Proof:} Let
\[
V(z)=\frac{1}{4}\sum_{i=1}^{n-1} (z_i^2-d_i^2)^2
\]
be a smooth non-negative function. We want to study the expression
taken by $\nabla V(z) f(x)$, where $f(z)$ is the vector field on
the right-hand side of (\ref{quantized.n.agent.system.z}). We
obtain:
\begin{eqnarray*}
&&\nabla_{z_i} V(z) \dot z_i   \\
&= & \left\{ \ba{ll}
z_1 (z_1^2-d_1^2) [-(k_1+1)\textrm{sgn}(z_1) \textrm{sgn}(|z_1|-d_1)\\+ k_2 \textrm{sgn}(z_2) \textrm{sgn}(|z_2|-d_2)] \\ \qquad i=1\\
\qquad \\
z_i(z_i^2-d_i^2)[-(k_i+1) \textrm{sgn}(z_i) \textrm{sgn}(|z_i|-d_{i}) \\ + \textrm{sgn}(z_{i-1}) \textrm{sgn}(|z_{i-1}|-d_{i-1})  \\
 +k_{i+1} \textrm{sgn}(z_{i+1}) \textrm{sgn}(|z_{i+1}|-d_{i+1})] \\ \qquad i=2,\ldots, n-2\\
 \qquad \\
z_{n-1}(z_{n-1}^2-d_{n-1}^2)[\textrm{sgn}(z_{n-2})
\textrm{sgn}(|z_{n-2}|-d_{n-2})
& \\
 -(k_{n-1}+1) \textrm{sgn}(z_{n-1}) \textrm{sgn}(|z_{n-1}|-d_{n-1})] \\
\qquad i=n-1 \ea \right.
\end{eqnarray*}

If $z\not \in {\cal D}$, i.e.\ if $z$ is not a point of
discontinuity for $f(z)$, then:
\begin{eqnarray*}
&&\nabla_{z_i} V(z) \dot z_i \\
&\le& \left\{ \ba{ll}
-(k_1+1-k_2) |z_1|\,|z_1^2-d_1^2|& i=1\\[2mm]
-(k_i-k_{i+1}) |z_i|\,|z_i^2-d_i^2| & i=2,\ldots, n-2\\[2mm]
-k_{n-1} |z_{n-1}|\,|z_{n-1}^2-d_{n-1}^2| & i=n-1 \ea \right.
\end{eqnarray*}
where we have exploited the fact that
$\textrm{sgn}(z_i^2-d_i^2)=\textrm{sgn}(|z_i|-d_i)$. Hence, if
(\ref{gains}) holds, then
\[
\nabla V(z) f(z)\le -\sum_{i=1}^{n-1}|z_i|\,|z_i^2-d_i^2|<0\;.
\]
If $z\in {\cal D}$, we look at the set
\[
\dot{\overline V}(z)=\{a\in \R: \exists v\in {\cal K}(f(z))\,\,\,
{\rm s.t.}\,\, \, a=\nabla V(z)\cdot v\}\;.
\]
We distinguish two cases, namely (i) $z\in {\cal E}\subseteq {\cal
D}$ and (ii) $z\in {\cal D}\setminus {\cal E}$. In case (i),
$\nabla V(z)=\mathbf{0}^T$, and therefore, $\dot{\overline
V}(z)=\{0\}$. In case (ii), there must exist at least one agent
such that $|z_i|\,|z_i^2-d_i^2|=0$ and at least one agent such
that $|z_j|\,|z_j^2-d_j^2|\ne 0$. Let ${\cal I}_1(z)$
(respectively, ${\cal I}_2(z)$) be the set of indices
corresponding to agents for which $|z_i|\,|z_i^2-d_i^2|=0$
($|z_j|\,|z_j^2-d_j^2|\ne 0$).
Clearly,  ${\cal I}_1(z)\cup {\cal I}_2(z)=\{1,2,\ldots,n-1\}$.\\
Since $\nabla_{z_i} V(z)=z_i(z_i^2-d_i^2)=0$ if $i\in {\cal
I}_1(z)$,   then
\[\ba{rcl}
\nabla V(z) \cdot v &=& \dst\sum_{i=1}^{n-1} z_i(z_i^2-d_i^2)
v_i\\
&=& \dst\sum_{i\in {\cal I}_2(z)} z_i(z_i^2-d_i^2) v_i\;. \ea
\]
Let $i\in {\cal I}_2(z)$ and $v\in {\cal K}(f(z))$.  In view of
(\ref{quantized.n.agent.system}), for $i=1,2,\ldots, n-1$, it
holds:
\begin{eqnarray*}
&&v_i\in \{\mu\in \R:
\mu=(2\lambda_1-1)\tilde k_{i+1}  \\
 && \qquad \qquad - (k_i+1) \textrm{sgn}(z_i) \textrm{sgn}(|z_i|-d_{i}),\,
\lambda_i\in [0,1]\}\;,
\end{eqnarray*}
with
\[
\tilde k_{i+1} =\left\{ \ba{ll}
k_{2} & i=1\\
1+k_{i+1} & i=2,\ldots, n-2\\
1 & i=n-1\;. \ea \right.
\]
Then
\[\ba{rcl}
\nabla V(z) \cdot v &=& \dst\sum_{i\in {\cal I}_2(z)}
z_i(z_i^2-d_i^2)
v_i\\
&\le& \dst\sum_{i\in {\cal I}_2(z)} -(k_i+1) |z_i| |z_i^2-d_{i}^2|
\\&& + \tilde k_{i+1} |z_i| |z_i^2-d_{i}^2| |2\lambda_j-1| \;. \ea
\]
By (\ref{gains}), $k_i+1-\tilde k_{i+1}\ge 1$ for all $i$, and
therefore, if $z\in {\cal D}\setminus {\cal E}$, then
\[\ba{rcl}
\nabla V(z) \cdot v &\le & -\dst\sum_{i\in {\cal I}_2(z)}
|z_i||z_i^2-d_i^2|<0 \;, \ea
\]
for all $v\in {\cal K}(f(z))$. This shows that for all $z\in {\cal
D}\setminus {\cal E}$, either $\max \dot{\overline V}(z)< 0$ or
$\dot{\overline V}(z)=\emptyset$. In summary,  for all $z\in
\R^{n-1}$, either $\max \dot{\overline V}(z)\le 0$ or
$\dot{\overline V}(z)=\emptyset$, and $0\in  \dot{\overline V}(z)$
if and only if $z\in {\cal E}$.

It is known (Lemma 1 in \cite{BaCe99}) that if $\varphi(t)$ is a
solution of the differential inclusion $\dot z\in {\cal K}(f(z))$,
then $\frac{d}{dt} V(\varphi(t))$ exists almost everywhere and
$\frac{d}{dt} V(\varphi(t))\in \dot{\overline V}(\varphi(t))$. We
conclude that $V(\varphi(t))$ is non-increasing. Let $z_0\in S$,
with $S\subset \R^{n-1}$  a compact and strongly invariant set for
(\ref{quantized.n.agent.system.z}). For any $z_0$, such a set
exists and includes the point $(d_1,d_2, \ldots, d_{n-1})\in {\cal
E}$ (hence $S\cap {\cal E}\ne \emptyset$), by definition of $V(z)$
and because $V(z)$ is non-increasing along the solutions of
(\ref{quantized.n.agent.system.z}). Since $\max\dot{\overline
V}(z)\le 0$ or $\dot{\overline V}(z)=\emptyset$ for all $z\in
\R^{n-1}$, then by the LaSalle invariance principle for
differential inclusions \cite{BaCe99, Co08}, any solution
$\varphi(t)$ to the differential inclusion converges to the
largest weakly invariant set in $S\cap \overline{{\cal E}}=S\cap
{\cal E}$  (${\cal E}$ is closed). Since the choice (\ref{gains})
yields that the gains $k_i$'s satisfy the condition in Lemma
\ref{lemma.equilibria}, ${\cal E}$ is the set of equilibria of
(\ref{quantized.n.agent.system.z}) (and therefore it is weakly
invariant) and since  $S\cap {\cal E}\ne \emptyset$, we conclude
that any solution $\varphi(t)$ converges to the set of points
$S\cap {\cal E}$. \hfill $\square$

Since the equilibrium set $\scr{E}$ contains those points for
which two agents coincide with each other, it is of interest to
characterize those initial conditions under which the asymptotic
positions of some of the agents become coincident. In the next
subsection, we use a three-agent formation as an example to show
how such analysis can be carried out.

\subsection{Trajectory based analysis}
\noindent We specialize the rigid formation examined before to the
case $n=3$. Letting $k_1=k_2=1$, the one-dimensional rigid
formation becomes:
\be\label{quantized.3.agent.system} \ba{rcr}
\dot x_1 &=& -\textrm{sgn}(x_1-x_2) \textrm{sgn}(|x_1-x_2|-d_1)\\
\dot x_2 &=& \textrm{sgn}(x_1-x_2) \textrm{sgn}(|x_1-x_2|-d_1)-\\
&& \textrm{sgn}(x_2-x_3) \textrm{sgn}(|x_2-x_3|-d_2)\\
\dot x_3 &=& \textrm{sgn}(x_2-x_3) \textrm{sgn}(|x_2-x_3|-d_2)\;.
\ea \ee

Let us express the system in the coordinates $z_1$, $z_2$, so as
to obtain:
\be\label{quantized.3.agent.system.z} \ba{rcr}
\dot z_1 &=& -2 \textrm{sgn}(z_1) \textrm{sgn}(|z_1|-d_1)+\textrm{sgn}(z_2) \textrm{sgn}(|z_2|-d_2)\\
\dot z_2 &=& \textrm{sgn}(z_1) \textrm{sgn}(|z_1|-d_1)-
2\textrm{sgn}(z_2) \textrm{sgn}(|z_2|-d_2)\;. \ea \ee
We study the solutions of the system above. In what follows, it is
useful to distinguish between two sets of points:
\begin{eqnarray*}
{\cal E}_1&=&\{z\in \R^2: |z_i|=d_i,\;i=1,2\}, \\ {\cal
E}_2&=&\{z\in \R^2: |z_i|=d_i\;\mbox{or}\; |z_i|=0,\;i=1,2\}\;.
\end{eqnarray*}
Clearly, ${\cal E}_1\subset {\cal E}_2$. We now prove that all the
solutions converge to the desired set ${\cal E}_1$ except for
solutions which originates on the $z_1$- or the $z_2$-axis:
\begin{theorem}\label{p1}
All Krasowskii solutions of (\ref{quantized.3.agent.system.z})
converge in finite time to the set ${\cal E}_2$. In particular,
the solutions which converge to the points $\{(d_1, 0), (0, d_2),
(-d_1, 0), (0, -d_2)\}$ must originate from the set of points
$\{z: z_1\cdot z_2=0, z\ne 0\}$. Moreover, the only solution which
converges to $(0,0)$ is the trivial solution which originates from
$(0,0)$.
\end{theorem}

\textit{Proof:} Because of the symmetry of $f(z)$, it suffices to
study the solutions which originate in the first quadrant only. In
the first quadrant we distinguish four regions: (i) ${\cal
R}_1=\{z\in \R^2: z_i\ge d_i, i=1,2\}$, (ii) ${\cal R}_2=\{z\in
\R^2: 0\le z_1<d_1,\; z_2\ge d_2\}$, (iii) ${\cal R}_3=\{z\in
\R^2: 0\le z_i<d_i,\;i=1,2\}$, (iv) ${\cal R}_4=\{z\in \R^2:
z_1\ge d_1,\; 0\le z_2< d_2\}$. Now we examine
the solutions originating in these regions.\\
(i) $z(0)\in {\cal R}_1$. If both $z_1(0)> d_1$ and $z_2(0)> d_2$,
then the system equations become
\[
\dot z_1 = -1\;,\quad \dot z_2 = -1
\]
and the solution satisfies $z_2(t)=z_1(t)+z_2(0)-z_1(0)$. In other
words, the solution evolves along the line of slop $+1$ and
intercept $z_2(0)-z_1(0)$. If $z_2(0)-z_1(0)=d_2-d_1$, then the
solution $z(t)$ converges to the point $z=(d_1,d_2)$ in finite
time. In particular $z(t_1)=(d_1,d_2)$ with
$t_1=z_1(0)-d_1=z_2(0)-d_2$. If $z_2(0)-z_1(0)>d_2-d_1$, then
$z(t)$ converges in finite time to the semi-axis $\{z:z_1=d_1,\;
z_2>d_2\}$. This is a set of points at which $f(z)$ is
discontinuous, since for $z_1\ge d_1$, $f(z)=(-1, -1)$,  and for
$z_1< d_1$, $f(z)=(3, -3)$. Since at these points $F(z)={\rm
co}\{(-1, -1), (3, -3)\}$,\footnote{Here ${\rm
co}\{v_1,\ldots,v_m\}$ denotes the smallest closed convex set
which contains $v_1,\ldots,v_m$.} and vectors in $F(z)$ intersect
the tangent space at the semi-axis in those points,  a sliding
mode along the semi-axis must occur. Since $\dot z(t)\in F(z(t))$,
we conclude that the sliding mode must satisfy the equations
\[
\dot z_1(t)=0,\; \dot z_2(t)=-\frac{3}{2}\;,
\]
and therefore, after a finite time, the solution converges to the
point $(d_1,d_2)$. On the other hand, if $z_2(0)-z_1(0)<d_2-d_1$,
then the solution reaches the ray $\{z:z_1>d_1,\; z_2=d_2\}$.
Similar considerations as before can show that a sliding mode
occurs along the ray and that it satisfies the equations
\[
\dot z_1(t)=-\frac{3}{2},\; \dot z_2(t)=0\;,
\]
and again convergence in finite time to $(d_1,d_2)$ is inferred.
Finally we examine the case $z(0)=(d_1,d_2)$. At the point
$(d_1,d_2)$,
\[
F(d_1,d_2)={\rm co}\{(-1, -1), (3, -3), (1, -1), (1,1)\}\;,
\]
i.e. $\mathbf{0}\in F(d_1,d_2)$ and $(d_1,d_2)$ is an equilibrium
point. Similarly as before, one shows that the solution which
originates from $(d_1,d_2)$ must stay in $(d_1,d_2)$.\\
(ii) $z(0)\in {\cal R}_2$. If $z_1(0)> 0$ and $z_2(0)>d_2$, then
the map $f(z)$ is equal to the vector $(3,-3)$ and the solution
$z(t)$ satisfies $z_2(t)=-z_1(t)+z_1(0)+z_2(0)$. If
$z_1(0)+z_2(0)=d_1+d_2$, then $z(t)$ converges to $(d_1,d_2)$,
while if $z_1(0)+z_2(0)=d_1+d_2$, it first converges to the
semi-axis $\{z:z_1=d_1,\; z_2>d_2\}$, and then it slides towards
$(d_1,d_2)$. When $z_1(0)+z_2(0)=d_1+d_2$, the solution reaches
the segment $\{z:0<z_1<d_1,\; z_2=d_2\}$. On this segment,
$F(z)={\rm co}\{(3, -3), (1, 1)\}$, and since this intersects the
tangent space at the segment, a sliding mode occurs. The sliding
mode obeys the equations
\[
\dot z_1(t)=\frac{3}{2},\; \dot z_2(t)=0\;,
\]
which show that the state reaches $(d_1,d_2)$.\\
If $z_1(0)=0$ and $z_2(0)>d_2$, then the initial condition lies on
another discontinuity surface of $f(z)$. Observe that, for those
points such that $-d_1<z_1\le 0 $ and $z_2> d_2$, $f(z)=(-1, -1)$.
Hence, $F(z)={\rm co}\{(-1, -1), (3, -3)\}$ intersects the tangent
space at the semi-axis in those points, and the solutions can
slide along the semi-axis until they reach the point $(0,d_2)$ and
stop, or can enter the region ${\cal R}_2\setminus
\{z:z_1=0,z_2>d_2\}$, and then converge to $(d_1,d_2)$, or they
can enter the region $\{z:
-d_1<z_1<0, z_2>d_2\}$ and converge to the point $(-d_1,d_2)$.\\
The point $(0,d_2)$ is an equilibrium, and if $z(0)=(0,d_2)$,
solutions stay at the equilibrium.\\
We review the remaining cases succinctly, as they are
qualitatively similar to the cases
examined above. \\
(iii) $z(0)\in {\cal R}_3$. If $z_i(0)> 0$ for $i=1,2$, then the
solutions converge to $(d_1,d_2)$ possibly sliding along the
segments $\{z:0<z_1\le d_1, z_2=d_2\}$ or $\{z:z_1=d_1, 0<z_2\le
d_2\}$. If $z_1(0)=0$ and $z_2(0)>0$, then the solution can
converge to the points $(-d_1, d_2)$, $(0, d_2)$ or $(d_1,d_2)$.
If $z_1(0)>0$ and $z_2(0)=0$, then the solutions can converge to
$(d_1,d_2)$, $(d_1,0)$ or $(d_1,-d_2)$. Finally, if $z_i(0)=0$ for
$i=1,2$, the solutions can converge to any of the points in ${\cal
E}_2$. In particular,
a possible solution is the one which remains in $(0,0)$.\\
(iv) $z(0)\in {\cal R}_4$. Solutions which start from initial
conditions  such that $z_1(0)>d_1$ and $z_2(0)>0$ converge to
$(d_1, d_2)$. If $z_1(0)=d_1$ and $z_2(0)>0$, then the solution
converge to $(d_1, d_2)$ possibly sliding on the segment
$\{z:z_1=d_1, 0<z_2<d_2\}$. If $z_1(0)>d_1$ and $z_2(0)=0$, the
solutions can converge to one of the three possible points: $(d_1,
-d_2)$, $(d_1, 0)$, $(d_1, d_2)$. \hfill $\square$

A few comments are in order:

\begin{itemize}
\item Sliding modes arise naturally for those situations in which, for
instance, the state reaches the semi-axis $\{z:z_1>d_1,
z_2=d_2\}$. This forces us to consider Krasowskii solutions rather
than Carath\'eodory solutions. On the other hand, the set of
Krasowskii solutions may be too large in some cases, as it is
evident for instance for those solutions which start on the $z_1$-
or $z_2$-axis.
\item The occurrence of sliding modes are not acceptable in practice as
they would require fast information transmission. A mechanism to
prevent sliding modes in the system
(\ref{quantized.3.agent.system}) can be introduced following
\cite{CePeFr10}. 
\end{itemize}

\section{Simulation results} \label{se:simulation}
\begin{figure}[ht]
 \centerline{\scalebox{.5}{\includegraphics{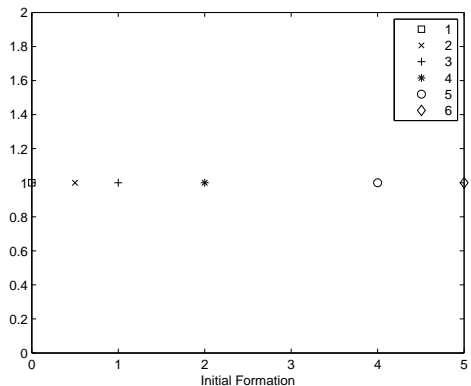}}}
    \caption{The initial shape of the 6-agent formation.}
    \label{fig1}
   \end{figure}

\begin{figure}[ht]
 \centerline{\scalebox{.5}{\includegraphics{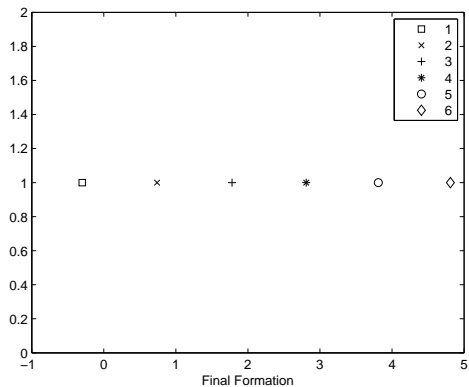}}}
    \caption{The final shape of the 6-agent formation.}
    \label{fig2}
   \end{figure}

\begin{figure}[ht]
 \centerline{\scalebox{.5}{\includegraphics{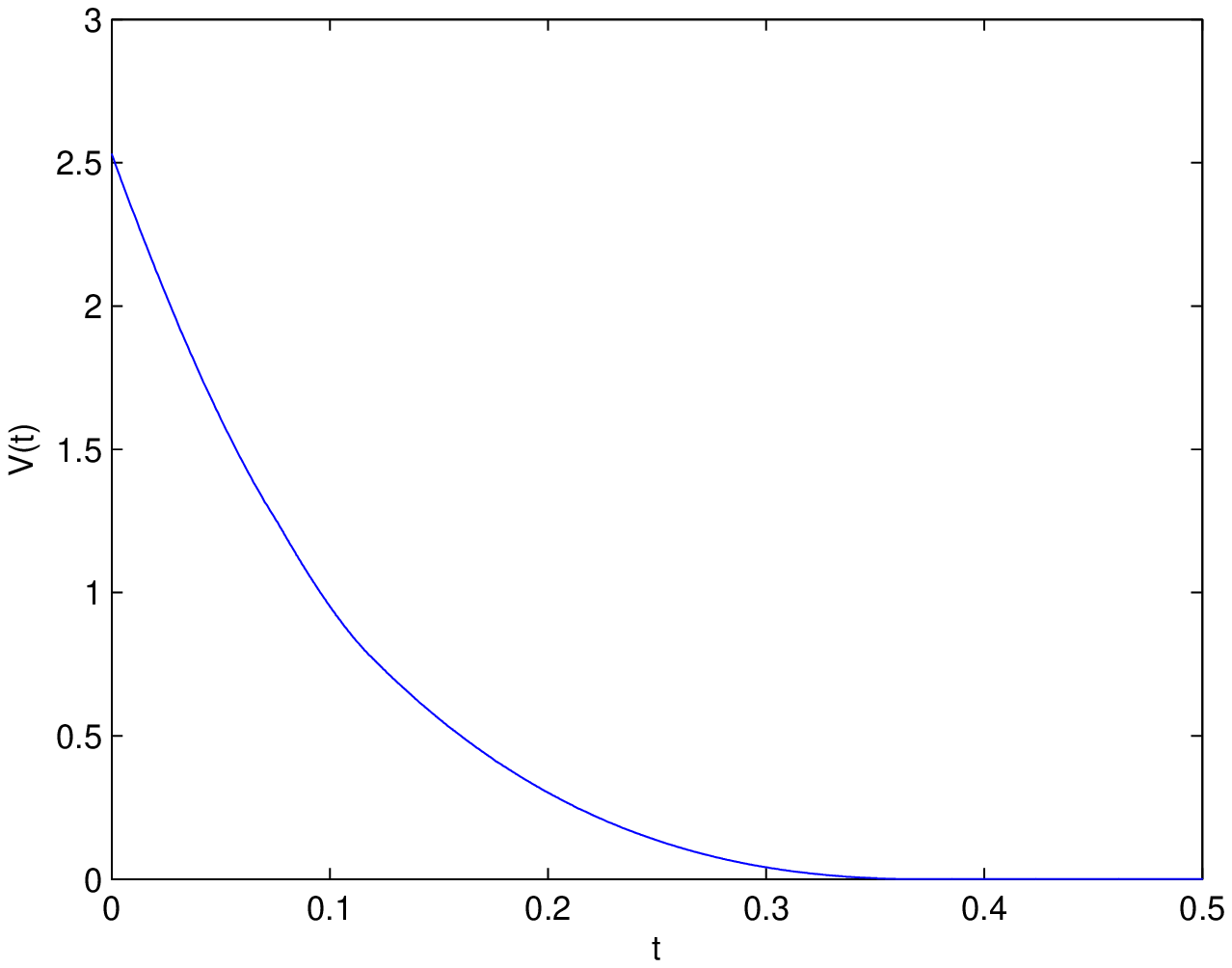}}}
    \caption{The curve of the Lyapunov function $V$.}
    \label{fig3}
   \end{figure}

\begin{figure}[ht]
 \centerline{\scalebox{.5}{\includegraphics{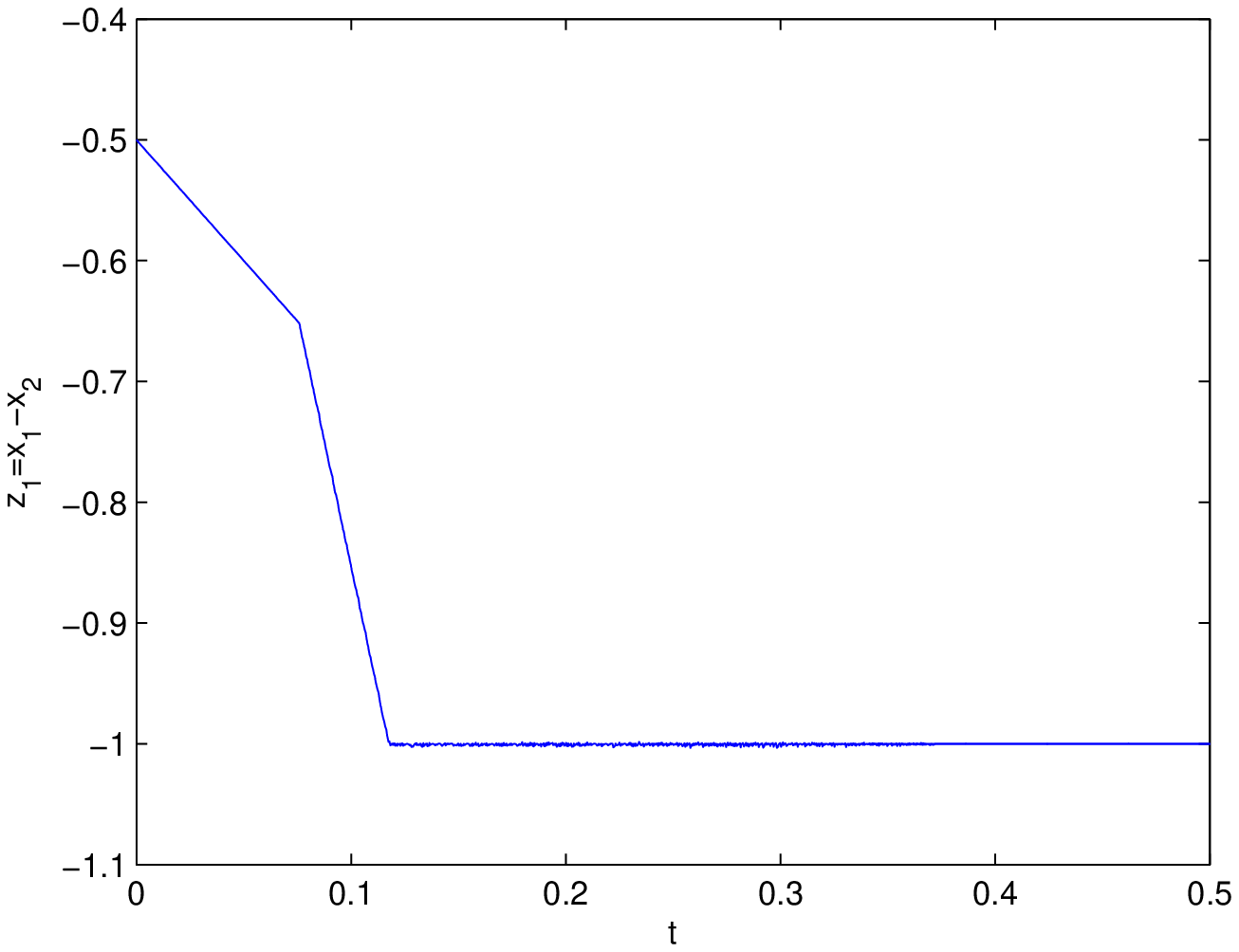}}}
    \caption{The dynamics of $z_1$.}
    \label{fig4}
   \end{figure}

In this section, we present simulation results for the guided
formation control with coarsely quantized information. We consider
a formation consisting of 6 agents, labeled by $1, \ldots, 6$. The
distance constraints are $|x_i-x_{i+1} |=1$, $i= 1, \ldots, 5$.
The initial positions of agents 1 to 6 are 0, 0.5, 1, 2, 4 and 5
respectively. Then the shape of the initial formation is shown in
figure \ref{fig1}.
We choose $k_1 =6$, $k_2 = 5$, $k_3=4$, $k_4=3$ and $k_5=2$ and
simulate the agents' motion under the control laws
(\ref{quantized.n.agent.system}). In figure \ref{fig2}, we show
the shape of the final formation.
To see how the shape evolves with time, we present the curve of
the Lyapunov function $V(z) = \frac{1}{4}\sum_{i=1}^5
(z_i^2-d_i^2)^2$ in figure \ref{fig3}.
Since our analysis has been carried out using Krasowskii
solutions, when we further look into the dynamics of $z$, it is
clear that the sliding mode may still happen when the Krasowskii
solution converges. But this effect due to the system's
non-smoothness  is within an acceptable level as shown in figure
\ref{fig4} which presents the curve of $z_1$.

\section{Concluding remarks} \label{se:conclusion}
In this paper, we have studied the problem of controlling a
one-dimensional guided formation using coarsely quantized
information. It has been shown that even when the guidance system
adopts quantizers that return only the one-bit sign information
about the quantized signal, the formation  can still converge to
the desired equilibrium under the proposed control law.

The point model we have used throughout the analysis is a
simplified description of vehicle dynamics. When more detailed
models are taken into consideration, we need to deal with
collision avoidance and other practical issues as well. So it is
of great interest to continue to study the same problem with more
sophisticated vehicle models and more physical constraints from
the applications.


\end{document}